\documentclass[12pt]{article}

\usepackage{amsmath}
\usepackage{amssymb}
\usepackage{vatola}
\usepackage{amsthm}
\allowdisplaybreaks[4]
\def\q{\quad}
\def\qq{\qquad}
\def\t{\hbox}
\def\qtq#1{\q\t{#1}\q}
\def\e{\equiv}
\def\f{\frac}
\def\b{\binom}
\def\ex{\text{ex}}
\let\pro=\proclaim
\let\endpro=\endproclaim

\begin{document}
\leftline{Czechoslovak Mathematical Journal 71(2021), no. 146, 351-372.}
\par\q\newline
 \centerline {\bf Ramsey numbers for trees II}
 \par\q
 \newline \centerline{Zhi-Hong Sun}\newline
\centerline{School of Mathematics
and Statistics}
\centerline{Huaiyin Normal University}
\centerline{Huaian, Jiangsu 223300, P.R. China} \centerline{Email:
zhsun@hytc.edu.cn} \centerline{URL:
http://maths.hytc.edu.cn/szh1.htm}
\par\q
\par {\it Abstract.} Let $r(G_1, G_2)$ be the Ramsey number of the two graphs
$G_1$ and $G_2$.  For $n_1\ge n_2\ge 1$ let $S(n_1,n_2)$ be the
double star given by
$V(S(n_1,n_2))=\{v_0,v_1,\ldots,v_{n_1},w_0,w_1,\ldots,w_{n_2}\}$
and $E(S(n_1,n_2))=\{v_0v_1,\ldots,v_0v_{n_1},v_0w_0,
w_0w_1,\ldots,w_0w_{n_2}\}$. In this paper we determine
$r(K_{1,m-1},$ $S(n_1,n_2))$ under certain conditions.
 For $n\ge 6$ let $T_n^3=S(n-5,3)$, $T_n^{''}=(V,E_2)$
 and $T_n^{'''} =(V,E_3)$, where $V=\{v_0,v_1,\ldots,v_{n-1}\}$,
$E_2=\{v_0v_1,\ldots,v_0v_{n-4},v_1v_{n-3},v_1v_{n-2},$
$v_2v_{n-1}\}$ and $E_3=\{v_0v_1,\ldots,$
$v_0v_{n-4},v_1v_{n-3},v_2v_{n-2},v_3v_{n-1}\}$. We also obtain
explicit formulas for $r$ $(K_{1,m-1},T_n)$, $r(T_m',T_n)$ $(n\ge
m+3)$, $r(T_n,T_n)$, $r(T_n',T_n)$ and $r(P_n,T_n)$, where
$T_n\in\{T_n'',T_n''',T_n^3\}$, $P_n$ is the path on $n$ vertices
and $T_n'$ is the unique tree with $n$ vertices and maximal degree
$n-2$.
\par\q
\par\q

 \let\thefootnote\relax\footnotetext {The author is supported by the
National Natural Science Foundation of China (Grant No. 11771173).}

\section*{1. Introduction}
 \par\q\  In this paper, all graphs are simple
graphs. For a graph $G=(V(G),E(G))$ let $e(G)=|E(G)|$ be the number
of edges in $G$, and let $\Delta(G)$ and $\delta(G)$ denote the
maximal degree and minimal degree of $G$, respectively.
\par For a graph $G$, as usual $\overline{G}$ denotes the
complement of $G$. Let $G\sb 1$ and $G\sb 2$ be two graphs. The
Ramsey number $r(G\sb 1, G\sb 2)$ is the smallest positive integer
$n$ such that, for every graph $G$ with $n$ vertices, either $G$
contains  a copy of $G\sb 1$
 or else $\overline{G}$ contains a copy of $G_2$.
 \par Let $\Bbb N$ be the set of positive integers. For
 $n\in\Bbb N$ with $n\ge 6$ let $T_n$ be a tree on $n$ vertices.
  As mentioned in [8],
 recently Zhao proved that $r(T_n,T_n)\le 2n-2$, which was conjectured by
 Burr and Erd\H os [1].
 \par Let $m,n\in\Bbb N$. For $n\ge 3$ let $K_{1,n-1}$ denote the unique tree on $n$
vertices with $\Delta(K_{1,n-1})=n-1$, and for $n\ge 4$ let
  $T_n'$ denote the unique tree on $n$ vertices
with $\Delta(T_n')=n-2$. In 1972, Harary [6] showed that for $m,n\ge
3$,
$$r(K_{1,m-1},K_{1,n-1})=\cases m+n-3&\t{if $2\nmid mn$,}
\\m+n-2&\t{if $2\mid mn$}.\endcases\tag 1.1$$
\par From [2, p.72], if
$G$ is a graph with $\delta(G)\ge n-1$, then $G$ contains every tree
on $n$ vertices. Using this fact, in 1995, Guo and Volkmann [5]
proved that for $n>m\ge 4$,
$$r(K_{1,m-1},T_n')=\cases m+n-3&\t{if $2\mid m(n-1)$,}
\\m+n-4&\t{if $2\nmid m(n-1)$.}
\endcases\tag 1.2$$
In 2012 the author [9] evaluated the Ramsey number $r(T_m,T_n^*)$
for $T_m\in\{P_m,K_{1,m-1},$ $T_m',T_m^*\}$, where $P_m$ is a path
on $m$ vertices and $T_n^*$ is the  tree on $n$ vertices with
 $V(T_n^*)=\{v_0,v_1,\ldots,v_{n-1}\}$ and
 $E(T_n^*)=\{v_0v_1,\ldots,v_0v_{n-3},v_{n-3}v_{n-2},$ $v_{n-2}v_{n-1}\}$.
 In particular, he proved that for $n>m\ge 7$,
    $$r(K_{1,m-1},T_n^*)=\cases m+n-3&\t{if $m-1\mid n-3$,}
\\m+n-4&\t{if $m-1\nmid n-3$.}\endcases\tag 1.3$$
For $n\ge 5$ let $T_n^1=(V,E_1)$ and $T_n^2=(V,E_2)$ be the trees on
$n$ vertices with
 $V=\{v_0,v_1,\ldots,v_{n-1}\}$,
 $E_1=\{v_0v_1,\ldots,v_0v_{n-3},v_{n-4}v_{n-2},v_{n-3}v_{n-1}\}$
and $E_2=\{v_0v_1,\ldots,v_0$ $v_{n-3},v_{n-3}v_{n-2},
v_{n-3}v_{n-1}\}$. Then $\Delta(T_n^1)=\Delta(T_n^2)$
$=\Delta(T_n^*)=n-3$. In [12], Sun, Wang and Wu proved that
$$r(K_{1,m-1},T_n^1)=r(K_{1,m-1},T_n^2)=m+n-4\ \t{for}\  n>m\ge 7 \ \t{and}
\ 2\mid mn.\tag 1.4$$

\par For $n_1,n_2\in\Bbb N$ with $n_1\ge n_2$, let $S(n_1,n_2)$ be the double star
given by
$$\align &V(S(n_1,n_2))=\{v_0,v_1,\ldots,v_{n_1},w_0,w_1,\ldots,w_{n_2}\},
\\&E(S(n_1,n_2))=\{v_0v_1,\ldots,v_0v_{n_1},v_0w_0, w_0w_1,\ldots,w_0w_{n_2}\}.
\endalign$$
We say that $v_0$ and $w_0$ are centers of $S(n_1,n_2)$. In [4],
Grossman, Harary and Klawe evaluated the Ramsey number
$r(S(n_1,n_2),S(n_1,n_2))$ under certain conditions. In particular,
they showed that for odd $n_1$ and $n_2=1,2$,
$$r(S(n_1,n_2),S(n_1,n_2))=max\{2n_1+1,n_1+2n_2+2\}.\tag 1.5$$ It is clear that $T_n'=S(n-3,1)$ and
$T_n^2=S(n-4,2)$. In this paper, we prove the following general
result:
$$\aligned &r(K_{1,m-1},S(n_1,n_2))
\\&=\cases m+n_1\qq\q\ \;\t{if $2\mid mn_1$, $n_1\ge m-2\ge n_2\ge 2$,
}\\\qq\qq\q\qq\t{and $n_1>m-5+n_2+\f{(n_2-1)(n_2-2)}{m-1-n_2}$,}
\\m-1+n_1\qq\t{if $2\nmid mn_1$, $n_1\ge m-2>n_2$,}\\\qq\qq\q
\qq\t{and $n_1>m-5+n_2+\f{(n_2-1)^2}{m-2-n_2}$.}
\endcases\endaligned\tag 1.6$$
Also,
$$r(K_{1,m-1},T_n^1)=m+n-5\qtq{for} n\ge m+2\ge 7 \q \t{and}\q 2\nmid
mn.\tag 1.7$$
 For $n\ge 6$ let $T_n^3=S(n-5,3)$, $T_n^{''}=(V,E_2)$ and $T_n^{'''}
=(V,E_3)$, where
$$\align &V=\{v_0,v_1,\ldots,v_{n-1}\},\q
 E_2=\{v_0v_1,\ldots,v_0v_{n-4},v_1v_{n-3},v_1v_{n-2},v_2v_{n-1}\},
\\&E_3=\{v_0v_1,\ldots,v_0v_{n-4},v_1v_{n-3},v_2v_{n-2},v_3v_{n-1}\}.
\endalign$$
Then $\Delta(T_n^3)=\Delta(T_n'')=\Delta(T_n''')=n-4$.
 In this paper, we evaluate  $r(K_{1,m-1},T_n)$ and
$r(T_m',T_n)$ for $T_n\in\{T_n'',T_n''',T_n^3\}$.  In particular, we
show that
$$\aligned& r(K_{1,m-1},T_n'')=r(K_{1,m-1},T_n''')
\\&=\cases m+n-5&\t{if $2\mid m(n-1)$, $m\ge 7$, $n\ge 15$
 and $n>m+1+\f 8{m-6}$,}
\\m+n-6&\t{if $2\nmid m(n-1)$ and $n\ge m+3\ge 9$,}
\endcases\endaligned\tag 1.8$$
and that for $m\ge 9$ and $n>m+2+\max\{0,\f {20-m}{m-8}\}$,
$$r(T_m',T_n'')=r(T_m',T_n''')=r(T_m',T_n^3)
=\cases m+n-5&\t{if $m-1\mid n-5$,}\\m+n-6&\t{if $m-1\nmid n-5$.}
\endcases\tag 1.9$$
 We also prove that for $m\ge 11$, $n\ge (m-3)^2+4$ and $m-1\nmid n-5$,
$$r(G_m,T_n)=m+n-6\ \t{for $G_m\in\{T_m^*,T_m^1,T_m^2\}$ and
$T_n\in\{T_n'',T_n''',T_n^3\}$}.\tag 1.10$$ In addition, we
establish the following results:
$$\align &r(T_n'',T_n'')=r(T_n'',T_n''')=r(T_n''',T_n''')=\cases
2n-9\q
 \t{if $2\mid n$ and $n>29$,}
 \\2n-8\q \t{if $2\nmid n$ and $n>22$},\endcases
\\&r(T_n^3,T_n'')=r(T_n^3,T_n''')=r(T_n^3,T_n^3)=2n-8\qtq{for} n>22,
\\&r(T_n'',T_n')=r(T_n''',T_n')=r(T_n^3,T_n')=2n-5\qtq{for} n\ge 10,
\\&r(T_n'',T_n^i)=r(T_n''',T_n^i)=r(T_n^3,T_n^i)=2n-7\q\t{for $n>16$
 and $i=1,2$,}
 \\&r(P_n,T_n'')=r(P_n,T_n''')=r(P_n,T_n^3)=2n-9\qtq{for} n\ge 33.
 \endalign$$
\par
In addition to the above notation, throughout this paper, we  use
the following notation: $[x]\f{\q}{\q}$the greatest integer not
exceeding $x$, $d(v)\f{\q}{\q}$the degree of the vertex $v$ in a
graph,  $d(u,v)\f{\q}{\q}$the distance between the two vertices $u$
and $v$ in a graph, $K_n\f{\q}{\q}$the complete graph on $n$
vertices, $G[V_1]\f{\q}{\q}$the subgraph of $G$ induced by vertices
in the set $V_1$, $G-V_1\f{\q}{\q}$the subgraph of $G$ obtained by
deleting vertices in $V_1$ and all edges incident with them,
$\Gamma(v)\f{\q}{\q}$the set of vertices adjacent to the vertex $v$,
 $e(V_1V_1')\f{\q}{\q}$the number of edges with one
endpoint in $V_1$ and the other endpoint in $V_1'$.

\section*{2. Basic lemmas}
\par\q
For a forbidden graph $L$ let $ex(p;L)$ be the maximal number of
edges in a graph of order $p$ not containing any copies of $L$. The
corresponding Tur\'an's problem is to evaluate $ex(p;L)$. Let
$p,n\in\Bbb N$ with $p\ge n\ge 2$. For a given tree $T_n$ on $n$
vertices,  it is difficult to determine the value of $ex(p;T_n)$.
The famous Erd\H os-S\'os
 conjecture asserts that
$ex(p;T_n)\le \f{(n-2)p}2$.
 Write
$p=k(n-1)+r$, where $k\in\Bbb N$ and $r\in\{0,1,\ldots,n-2\}$.  In
1975 Faudree and Schelp [3] showed that
$$ex(p;P_n)=k\binom {n-1}2+\binom r2=\f{(n-2)p-r(n-1-r)}2.\tag 2.1$$
In [10-12], the author and his coauthors determined $\ex(p;T_n)$ for
$T_n\in\{T_n',T_n^*,T_n^1,$ $T_n^2,T_n^3,T_n'',$ $T_n'''\}$.

\pro{Lemma 2.1 ([9, Lemma 2.1])} Let $G_1$ and $G_2$ be two graphs.
Suppose that $p\in\Bbb N, p\ge max\{|V(G_1)|$, $|V(G_2)|\}$ and
$ex(p;G_1)+ex(p;G_2)<\b p2.$ Then $r(G_1,G_2)\le p.$\endpro
 Proof. Let $G$ be a graph of order $p$. If
$e(G)\le ex(p;G_1)$ and $e(\overline{G})\le ex(p;G_2)$, then
$ex(p;G_1)+ex(p;G_2)\ge e(G)+e(\overline{G})=\b p2.$ This
contradicts the assumption. Hence, either $e(G)>ex(p;G_1)$ or
$e(\overline{G})>ex(p;G_2).$ Therefore, $G$ contains a copy of $G_1$
or $\overline{G}$ contains a copy of $G_2$. This shows that
$r(G_1,G_2)\le |V(G)|=p.$ So the lemma is proved.
\par\q
\pro{Lemma 2.2} Let $k,p\in\Bbb N$ with $p\ge k+1$. Then there
exists a $k-$regular graph of order $p$ if and only if $2\mid kp$.
\endpro
\par This is a known result. See for example [11, Corollary 2.1].
\par\q
 \pro{Lemma 2.3 ([9, Lemma 2.3])} Let $G_1$ and $G_2$ be two
graphs with $\Delta(G_1)=d_1\ge 2$ and $\Delta(G_2)=d_2\ge 2$. Then
\par $(\t{\rm i})$ $r(G_1,G_2)\ge d_1+d_2-(1-(-1)^{(d_1-1)(d_2-1)})/2$.
\par $(\t{\rm ii})$ Suppose that $G_1$ is a connected graph of order
$m$ and $d_1<d_2\le m$. Then $r(G_1,G_2)\ge 2d_2-1\ge d_1+d_2$.
\par $(\t{\rm iii})$ Suppose that $G_1$ is a connected graph of order
$m$ and $d_2>m$.  If one of the conditions:
\par \qq$(\t{\rm 1})$ $2\mid (d_1+d_2-m)$,
\par \qq$(\t{\rm 2})$ $d_1\not=m-1$,
\par \qq$(\t{\rm 3})$ $G_2$ has two
vertices $u$ and $v$ such that $d(v)=\Delta(G_2)$ and
\par\qq \q\  $d(u,v)=3$
\newline holds, then $r(G_1,G_2)\ge d_1+d_2$.
\endpro
\par\q
 \pro{Lemma 2.4} Let $p,n\in\Bbb N$ with $p\geq n-1\geq 1$. Then
$ex(p;K_{1,n-1})=[\f{(n-2)p}2]$.
\endpro
\par This is a known result. See for example [11, Theorem 2.1].
\par\q
\pro{Lemma 2.5 ([11, Theorem 3.1])} Let $p,n\in\Bbb N$ with $p\geq
n\geq 5$, and let $r\in\{ 0,1,\ldots,n-2\}$ be given by $p\e
r\mod{n-1}$. Then
$$ex(p;T_n')= \cases [\f{(n-2)(p-1)-r-1}2] &\t{if}\
n\ge 7\ \t{and}\
2\le r\le n-4,\\
\f{(n-2)p-r(n-1-r)}2&\t{otherwise.}\endcases$$
\endpro

\pro{Lemma 2.6 ([12, Theorems 2.1 and 3.1])} Suppose that
$p,n\in\Bbb N,\ p\geq n-1\geq 4$ and $p=k(n-1)+r$, where $k\in \Bbb
N$ and $r\in\{ 0,1,\ldots,n-2\}$. For $i=1$ or $2$,
 $$\aligned &\ex(p;T_n^i)\\&=
\max\Big\{\Big[\f{(n-2)p}2\Big]-(n-1+r),\f{(n-2)p-r(n-1-r)}2\Big\}
\\&= \cases [\f{(n-2)p}2]-(n-1+r) &\t{if}\ n\ge 16\ \t{and}\
3\le r\le n-6, \ \t{or if}\\&\ \ 13\le n\le 15 \ \t{and}\ 4\le r\le n-7,\\
\f{(n-2)p-r(n-1-r)}2&\t{otherwise.}\endcases\endaligned$$
\endpro

\pro{Lemma 2.7 ([10, Theorems 3.1 and 5.1])} Let $p,n \in \Bbb N,
p\geq n\ge 10$, $p=k(n-1)+r$, $k\in\Bbb N$ and
$r\in\{0,1,\ldots,n-2\}$. Then
$$\align ex(p;T_n'')=ex(p;T_n''')&=
\f{(n-2)p-r(n-1-r)}2\\&\q+\max\big\{0,\big[\f{r(n-4-r)
-3(n-1)}2\big]\big\}.\endalign$$\endpro

\pro{Lemma 2.8 ([10, Lemmas 4.6 and 4.7])} Let $n \in \Bbb N$ with
$n\ge 15$. Then
$$\ex(2n-9;T_n^3)=n^2-10n+24+\max\Big\{\Big[\f n2\Big],13\Big\}$$
and
$$\ex(2n-8;T_n^3)=n^2-9n+29+\max\Big\{\Big[\f {n-37}4\Big],0\Big\}.$$
\endpro

\pro{Lemma 2.9 ([10, Theorems 4.1-4.5])} Let $p,n \in \Bbb N, p\ge
n\ge 10$, $p=k(n-1)+r$, $k\in\Bbb N$ and $r\in\{0,1,\ldots,n-2\}$.
\par $(\t{\rm i})$ If $r\in\{0,1,2,n-6,n-5,n-4,n-3,n-2\}$, then
$$\ex(p;T_n^3)=\f{(n-2)p-r(n-1-r)}2.$$
\par $(\t{\rm ii})$ If $n\ge 15$ and $r\in\{3,4,\ldots,n-9\}$, then
$$\align &ex(p;T_n^3)\\&=
\f{(n-2)p-r(n-1-r)}2+\max\Big\{0,\Big[\f{r(n-4-r)
-3(n-1)}2\Big]\Big\}.\endalign$$
\par $(\t{\rm iii})$ If $n\ge 15$ and $r=n-8$, then
$$ex(p;T_n^3)=\f{(n-2)p-7n+30}2+\max\Big\{\Big[\f
n2\Big],13\Big\}.$$
\par $(\t{\rm iv})$ If $n\ge 15$ and $r=n-7$, then
$$ex(p;T_n^3)=\f{(n-2)p-6(n-7)}2+\max\Big\{\Big[\f
{n-37}4\Big],0\Big\}.$$
\endpro
\pro{Lemma 2.10} Let $n\in\Bbb N$, $n\ge 10$ and
$T_n\in\{T_n'',T_n''',T_n^3\}$. Assume that $p=k(n-1)+r$ with
$k\in\Bbb N$ and $r\in\{0,1,\ldots,n-2\}$. Then
$$\ex(p;T_n)\le \f{(n-2)p}2-\min\Big\{n-1+r,\f{r(n-1-r)}2\Big\}.$$
\endpro
Proof. This is immediate from [10, Lemmas 2.8, 3.1, 4.1 and 5.1].
\par\q
\pro{Lemma 2.11 ([11, Theorems 4.1-4.3])} Let $p,n\in\Bbb N$, $p\ge
n\ge 6$ and  $p=k(n-1)+r$ with $k\in\Bbb N$ and
$r\in\{0,1,n-5,n-4,n-3,n-2\}.$ Then $$\aligned ex(p;T_n^*)&=\cases
\f{(n-2)(p-2)}2+1&\t{if $n>6$ and $r=n-5$,}
\\\f{(n-2)p-r(n-1-r)}2&\t{otherwise.}
\endcases\endaligned$$\endpro
\pro{Lemma 2.12 ([11, Theorem 4.4])} Let $p,n\in\Bbb N$, $p\ge n\ge
11$, $r\in\{2,3,\ldots,n-6\}$ and $p\e r\mod{n-1}$. Let
$t\in\{0,1,\ldots,r+1\}$ be given by $n-3\e t\mod{r+2}$. Then
$$\align &\ex(p;T_n^*)\\&=\cases [\f{(n-2)(p-1)-2r-t-3}2
]&\t{if $r\ge 4$ and $2\le t\le r-1$,}\\\f{(n-2)(p-1)-t(r+2-t)-r-1}2
&\t{otherwise}.\endcases\endalign$$\endpro

\section*{3. Formulas for $r(T_n,T_n''),r(T_n,T_n''')$ and
$r(T_n,T_n^3)$}
\par\q
\pro{Theorem 3.1} Let $n\in\Bbb N$. Then
$$ r(T_n'',T_n'')=r(T_n'',T_n''')=r(T_n''',T_n''')=\cases
2n-9\q
 \t{if $2\mid n$ and $n>29$,}
 \\2n-8\q \t{if $2\nmid n$ and $n>22$}.\endcases$$\endpro
Proof. Suppose that $T_n,T_n^0\in \{T_n'',T_n'''\}$. By Lemma 2.7,
$$\align \ex(2n-9;T_n)&=\f{(2n-9)(n-5)-(n-29)}2+\max\Big\{0,
\Big[\f{n-29}2\Big]\Big\}
\\&=\Big[\f{(2n-9)(n-5)}2\Big]\q\t{for $n\ge 29$}.
\endalign$$
Hence, for $n\in\{30,32,34,\ldots\}$,
$$\align &\ex(2n-9;T_n)+\ex(2n-9;T_n^0)
\\&=2\Big[\f{(2n-9)(n-5)}2\Big]=(2n-9)(n-5)-1<\b{2n-9}2.
\endalign$$
Applying Lemma 2.1 yields $r(T_n,T_n^0)\le 2n-9$. On the other
hand, appealing to Lemma 2.3(i),
$$r(T_n,T_n^0)\ge n-4+n-4-\f{1-(-1)^{(n-5)(n-5)}}2=2n-9.$$ Therefore
$r(T_n,T_n^0)=2n-9$ for $n\in\{30,32,34,\ldots\}$.

\par Now assume that $2\nmid n$ and $n>22$. By Lemma 2.7,
$$\ex(2n-8;T_n)=\f{(n-2)(2n-8)-6(n-7)}2=n^2-9n+29.$$
Thus,
$$\align&\ex(2n-8;T_n)+\ex(2n-8;T_n^0)
\\&=2(n^2-9n+29)<2n^2-17n+36=\b{2n-8}2.\endalign$$
 Hence $r(T_n,T_n^0)\le 2n-8$
by  Lemma 2.1. By Lemma 2.2, we may construct a regular graph $G$ of
order $2n-9$ with degree $n-5.$ Clearly $\overline{G}$ is also a
regular graph with degree $n-5.$ Since
$\Delta(T_n)=\Delta(T_n^0)=n-4$,  both $G$ and $\overline{G}$ do not
contain any copies of $T_n$ and $T_n^0$. Therefore,
$r(T_n,T_n^0)>2n-9$ and so $r(T_n,T_n^0)=2n-8$. This completes the
proof.
\par\q
\pro{Theorem 3.2} Let $n\in\Bbb N$ with $n>22$. Then
$$r(T_n^3,T_n'')=r(T_n^3,T_n''')=r(T_n^3,T_n^3)=2n-8.$$
\endpro
Proof. Let $T_n\in\{T_n'',T_n''',T_n^3\}$. When $n$ is odd, using
Lemma 2.3(i) we see that $r(T_n^3,T_n) \ge n-4+n-4=2n-8$.   When $n$
is even, we may construct a regular graph $H$ with degree $n-10$ and
$V(H)=\{v_1,\ldots,v_{n-6}\}$. Let $G_0$ be a graph given by
$$V(G_0)=\{v_0,v_1,\ldots,v_{n-4},u_1,\ldots,u_{n-6}\}$$ and
$$\align &E(G_0)\\&=E(H)\cup \{v_0v_1,\ldots,v_0v_{n-4},v_1v_{n-5},
\ldots, v_{n-6}v_{n-5}, v_1v_{n-4},\ldots,v_{n-5}v_{n-4},
\\&\q v_1u_1,v_1u_2,v_2u_1,v_2u_2,\ldots,v_{n-7}u_{n-7},v_{n-7}u_{n-6},
v_{n-6}u_{n-7},v_{n-6}u_{n-6},
\\&\q u_1u_2,\ldots, u_1u_{n-6},u_2u_3,\ldots,u_2u_{n-6},
u_3u_{n-6},\ldots,u_{n-7}u_{n-6}\big\}.
\endalign$$
Then $d(v_0)=d(v_{n-5})=d(v_{n-4})=n-4$ and
$d(v_1)=\cdots=d(v_{n-6})=d(u_1)=\cdots=d(u_{n-6})=n-5$. Clearly
$|V(G_0)|=2n-9$ and $G_0$ does not contain any copies of $T_n^3$.
Since $\Delta(\overline{G_0})=n-5$ and $\Delta(T_n)=n-4$,
$\overline{G_0}$ does not contain any copies of $T_n$. Thus,
$r(T_n^3,T_n)\ge |V(G_0)|+1=2n-8$.
\par From Lemma 2.7, $\ex(2n-8;T_n'')=\ex(2n-8;T_n''')
=n^2-9n+29$. By Lemma 2.8,
$\ex(2n-8;T_n^3)=n^2-9n+29+\max\big\{0,[\f{n-37}4]\big\}$. Thus,
$$\align &\ex(2n-8;T_n^3)+\ex(2n-8;T_n)\\&
\le 2n^2-18n+58+2\max\{0,[\f{n-37}4]\}
\\&<2n^2-18n+58+n-22=\b{2n-8}2.\endalign$$
Hence, applying Lemma 2.1 gives $r(T_n^3,T_n)\le 2n-8$ and so
$r(T_n^3,T_n)=2n-8$ as claimed.
\par\q
 \pro{Theorem 3.3} Let
$n\in\Bbb N$ with $n\ge 10$. Then
$$r(T_n'',T_n')=r(T_n''',T_n')=r(T_n^3,T_n')=2n-5.$$
\endpro
Proof. Let $T_n\in\{T_n'',T_n''',T_n^3\}$. Since $\Delta(T_n)=n-4$
and $\Delta(T_n')=n-2$, using Lemma 2.3(ii) we see that
$r(T_n,T_n')\ge 2(n-2)-1=2n-5$. By Lemmas 2.5, 2.7 and 2.9,
$$\align &\ex(2n-5;T_n)+\ex(2n-5;T_n')
\\&=\f{(n-2)(2n-5)-3(n-4)}2+ \big[\f{(n-2)(2n-6)-(n-3)}2\big]
\\&=\big[\f{4n^2-23n+37}2\big] <\f{4n^2-22n+30}2=\b{2n-5}2.
\endalign$$
Hence, $r(T_n,T_n')\le 2n-5$ by Lemma 2.1. Therefore,
$r(T_n,T_n')=2n-5$ as claimed.
\par\q
\pro{Theorem 3.4} Let $n\in\Bbb N$, $n>16$ and $i\in\{1,2\}$. Then
$$r(T_n'',T_n^i)=r(T_n''',T_n^i)=r(T_n^3,T_n^i)=2n-7.$$
\endpro
Proof. Let $T_n\in\{T_n'',T_n''',T_n^3\}$. Since $\Delta(T_n)=n-4$
and $\Delta(T_n^i)=n-3$, using Lemma 2.3(ii) we see that
$r(T_n,T_n^i)\ge 2(n-3)-1=2n-7$. From Lemmas 2.6, 2.7 and 2.9,
$$\align &\ex(2n-7;T_n)+\ex(2n-7;T_n^i)
\\&=\f{(n-2)(2n-7)-5(n-6)}2+ \Big[\f{(n-2)(2n-7)}2\Big]-(2n-7)
 \\&=\Big[\f{4n^2-31n+72}2\Big]<\f{4n^2-30n+56}2=\b{2n-7}2.
\endalign$$
Hence, $r(T_n,T_n^i)\le 2n-7$ by Lemma 2.1. Therefore,
$r(T_n,T_n^i)=2n-7$ as claimed.
\par\q
\pro{Theorem 3.5} Let $n\in\Bbb N$ with $n\ge 10$. Then
$r(T_n,T_n^*)=2n-5$ for $T_n\in\{T_n'',T_n''',T_n^3\}$.
\endpro
Proof. By Lemmas 2.7 and 2.9,
$\ex(2n-5;T_n)=\f{(n-2)(2n-5)-3(n-4)}2=n^2-6n+11<n^2-5n+4$. Thus the
result follows from [9, Lemma 3.1].
\par\q
\par{\bf Remark 3.1} By [9, Theorem 6.3 (with $m=n$ and $a=2$)],
$r(T_n,$ $K_{1,n-1})=2n-3$ for $n\ge 6$ and
$T_n\in\{T_n'',T_n''',T_n^3\}$.
\par\q
\pro{Theorem 3.6} Let $n\in\Bbb N$.  Then
$r(P_n,T_n'')=r(P_n,T_n''')=2n-9$ for $n\ge 30$ and
$r(P_n,T_n^3)=2n-9$ for $n\ge 33$.
\endpro
Proof. Suppose that $n\ge 30$ and $T_n\in\{T_n'',T_n''',T_n^3\}$.
Since $\Delta(T_n)=n-4$ and $\Delta(P_n)=2$, appealing to Lemma
2.3(ii) we obtain $r(P_n,T_n)\ge 2(n-4)-1=2n-9$. By (2.1) and Lemma
2.7, for $T_n\in\{T_n'',T_n'''\}$,
$$\align &\ex(2n-9;P_n)+\ex(2n-9;T_n)
\\&=\f{(n-2)(2n-9)-7(n-8)}2+\Big[\f{(2n-9)(n-5)}2\Big]
=\Big[\f{4n^2-39n+119}2\Big]
\\&<\f{4n^2-38n+90}2=\b{2n-9}2.\endalign$$
Hence, applying Lemma 2.1 gives $r(P_n,T_n)\le 2n-9$ and so
$r(P_n,T_n)=2n-9$.
\par Now assume that $n\ge 33$. From (2.1) and Lemma 2.8,
$$\align &\ex(2n-9;P_n)+\ex(2n-9;T_n^3)
\\&=\f{(n-2)(2n-9)-7(n-8)}2+n^2-10n+24+\Big[\f n2\Big]
\\&=2n^2-20n+61+\big[\f n2\big] <2n^2-19n+45=\b{2n-9}2.\endalign$$
Hence, $r(P_n,T_n^3)\le 2n-9$ by Lemma 2.1 and so
$r(P_n,T_n^3)=2n-9$.

\section*{4. Formulas for $r(K_{1,m-1},S(n_1,n_2))$,
$r(K_{1,m-1},$\\$T_n^1)$, $r(K_{1,m-1},T_n^{''})$ and
$r(K_{1,m-1},T_n^{'''})$} \pro{Theorem 4.1} Let $m,n_1,n_2\in\Bbb N$
with $n_1\ge m-2\ge n_2\ge 2$ and $2\mid mn_1$. If
$n_1>m-5+n_2+\f{(n_2-1)(n_2-2)}{m-1-n_2}$, then
$r(K_{1,m-1},S(n_1,n_2))=m+n_1.$
\endpro
Proof. Since $\Delta(S(n_1,n_2))=n_1+1$, from Lemma 2.3(i) we see
that
$$r(K_{1,m-1},S(n_1,n_2))\ge
m-1+n_1+1-\f{1-(-1)^{(m-2)n_1}}2=m+n_1.$$ Now we show that
$r(K_{1,m-1},S(n_1,n_2))\le m+n_1$. Let $G$ be a graph of order
$m+n_1$ such that $\overline{G}$ does not contain any copies of
$K_{1,m-1}$. That is, $\Delta(\overline{G})\le m-2$. We show that
$G$ contains a copy of $S(n_1,n_2)$. Clearly
$$\delta(G)=m+n_1-1-\Delta(\overline{G})\ge m+n_1-1-(m-2)=n_1+1.$$
Suppose that $\Delta(G)=n_1+1+s$,
 $v_0\in V(G)$, $d(v_0)=\Delta(G)$,
$\Gamma(v_0)=\{v_1,\ldots,v_{n_1+1+s}\}$,
$V_1=\{v_0\}\cup\Gamma(v_0)$ and $V_1'=V(G)-V_1$. Then
$|V_1'|=m-2-s$. For $i=1,2,\ldots,n_1+1+s$, we have
$$|V_1'|+1+|\Gamma(v_i)\cap\Gamma(v_0)|\ge d(v_i)\ge \delta(G)\ge
n_1+1$$ and so
$$|\Gamma(v_i)\cap\Gamma(v_0)|\ge
n_1-|V_1'|=n_1-(m-2)+s\ge s.$$ For $s\ge n_2$ we have
$|\Gamma(v_i)\cap\Gamma(v_0)|\ge s\ge n_2$ and
$|\Gamma(v_0)|-n_2=n_1+1+s-n_2\ge n_1+1$. Hence $G[V_1]$ contains a
copy of $S(n_1,n_2)$ with centers $v_0$ and $v_i$.
\par Now assume that $s<n_2$ and $V_1'=V(G)-V_1=\{u_1,\ldots,u_{m-2-s}\}$.
It is clear that for $i=1,2,\ldots,m-2-s$,
$$m-3-s+|\Gamma(u_i)\cap
\Gamma(v_0)|=|V_1'|-1+|\Gamma(u_i)\cap \Gamma(v_0)|\ge d(u_i)\ge
\delta(G)\ge n_1+1$$ and so
 $|\Gamma(u_i)\cap
\Gamma(v_0)|\ge n_1-(m-4-s)$. It then follows that $e(V_1V_1')\ge
(m-2-s)(n_1-(m-4-s))$. By the assumption,
$$n_1>m-5+n_2+\f{(n_2-2)(n_2-1)}{m-1-n_2}
\ge m-5+n_2-2s+\f{(n_2-2)(n_2-1-s)}{m-1-n_2}.$$ Thus, $(m-1-n_2)n_1>
(m-2-s)(m-4-s)+(s+1)(n_2-s-1)$ and so $e(V_1V_1')\ge
(m-2-s)(n_1-(m-4-s))>(n_1+1+s)(n_2-s-1)$. Therefore
$|\Gamma(v_i)\cap V_1'|\ge n_2-s$ for some $v_i\in \Gamma(v_0)$.
From the above, $|\Gamma(v_i)\cap\Gamma(v_0)|\ge s$. Thus, $G$
contains a copy of $S(n_1,n_2)$ with centers $v_0$ and $v_i$.
Therefore $r(K_{1,m-1},S(n_1,n_2))\le m+n_1$ and so the theorem is
proved. \par\q
\pro{Corollary 4.1} Let $m,n\in\Bbb N$, $n-2\ge m\ge
4$ and $2\mid mn$. Then $r(K_{1,m-1},T_n^2)=m+n-4$.
\endpro
Proof. Since $T_n^2=S(n-4,2)$, putting $n_1=n-4$ and $n_2=2$ in
Theorem 4.1 yields the result.
\par\q
\pro{Corollary 4.2} Let $m,n\in\Bbb N$, $m\ge 5$, $n>m+3+\f 2{m-4}$
and $2\mid m(n-1)$. Then $r(K_{1,m-1},T_n^3)=m+n-5$.
\endpro
Proof. Since $T_n^3=S(n-5,3)$, taking $n_1=n-5$ and $n_2=3$ in
Theorem 4.1 deduces the result.
\par\q
\pro{Theorem 4.2} Let $m,n_1,n_2\in\Bbb N$, $n_1\ge m-2>n_2$ and
$2\nmid mn_1$. If $n_1>m-5+n_2+\f{(n_2-1)^2}{m-2-n_2}$, then
$r(K_{1,m-1},S(n_1,n_2))=m-1+n_1.$
\endpro
Proof. Since $\Delta(S(n_1,n_2))=n_1+1$, from Lemma 2.3(i) we see
that
$$r(K_{1,m-1},S(n_1,n_2))\ge
m-1+n_1+1-\f{1-(-1)^{(m-2)n_1}}2=m-1+n_1.$$ Now we show that
$r(K_{1,m-1},S(n_1,n_2))\le m-1+n_1$. Let $G$ be a graph of order
$m-1+n_1$ such that $\overline{G}$ does not contain any copies of
$K_{1,m-1}$. We need to show that $G$ contains a copy of
$S(n_1,n_2)$. Clearly $\Delta(\overline{G})\le m-2$ and so
$\delta(G)=m-2+n_1-\Delta(\overline{G})\ge n_1$. Since $2\nmid
mn_1$, there is no regular graph of order $m-1+n_1$ with degree
$n_1$ by Euler's theorem. Hence $\Delta(G)\ge \delta(G)+1\ge n_1+1$.
Suppose that $\Delta(G)=n_1+1+s$,
 $v_0\in V(G)$, $d(v_0)=\Delta(G)$,
$\Gamma(v_0)=\{v_1,\ldots,v_{n_1+1+s}\}$,
$V_1=\{v_0\}\cup\Gamma(v_0)$ and $V_1'=V(G)-V_1$. Then
$|V_1'|=m-3-s$. For $v_i\in\Gamma(v_0)$,  $d(v_i)\ge \delta(G)\ge
n_1$ and so $|\Gamma(v_i)\cap\Gamma(v_0)|+1+|V_1'|\ge d(v_i)\ge
n_1$. Thus, $$|\Gamma(v_i)\cap\Gamma(v_0)|\ge
n_1-1-|V_1'|=n_1-(m-2)+s\ge s.$$ Hence $G[V_1]$ contains a copy of
$S(n_1,n_2)$ with centers $v_0$ and $v_i$ for $s\ge n_2$.
\par Now assume that $s<n_2$ and $V_1'=\{u_1,\ldots,u_{m-3-s}\}$.
As $d(u_i)\ge n_1$, we see that $|\Gamma(u_i)\cap \Gamma(v_0)|\ge
n_1-(m-4-s)$ and so $e(V_1V_1')\ge (m-3-s)(n_1-(m-4-s))$. Since
$$n_1>m-5+n_2+\f{(n_2-1)^2}{m-2-n_2}\ge
m-5+n_2-2s+\f{(n_2-1)(n_2-1-s)}{m-2-n_2},$$ we get $(m-2-n_2)n_1>
(m-3-s)(m-4-s)+(s+1)(n_2-s-1)$. Hence
$$\align e(V_1V_1')&\ge
(m-3-s)(n_1-(m-4-s))\\&>(m-3-s)n_1-(m-2-n_2)n_1+(s+1)(n_2-s-1)
\\&=(n_1+1+s)(n_2-s-1).\endalign$$
 Therefore $|\Gamma(v_i)\cap V_1'|\ge n_2-s$
for some $v_i\in \Gamma(v_0)$. From the above,
$|\Gamma(v_i)\cap\Gamma(v_0)|\ge s$. Thus, $G$ contains a copy of
$S(n_1,n_2)$ with centers $v_0$ and $v_i$. Therefore,
$r(K_{1,m-1},S(n_1,n_2))$ $\le m-1+n_1$ and so the theorem is
proved.
\par\q
\pro{Corollary 4.3} Let $m,n\in\Bbb N$, $m\ge 5$, $n> m+1+\f 1{m-4}$
and $2\nmid mn$. Then $r(K_{1,m-1},T_n^2)=m+n-5$.
\endpro
Proof. Since $T_n^2=S(n-4,2)$, putting $n_1=n-4$ and $n_2=2$ in
Theorem 4.2 yields the result.
\par\q

\pro{Corollary 4.4} Let $m,n\in\Bbb N$, $m\ge 6$, $n>m+3+\f 4{m-5}$
and $2\nmid m(n-1)$. Then $r(K_{1,m-1},T_n^3)=m+n-6$.
\endpro
Proof. Since $T_n^3=S(n-5,3)$, putting $n_1=n-5$ and $n_2=3$ in
Theorem 4.2 we deduce the result.
\par\q
\pro{Theorem 4.3} Let $m,n\in\Bbb N$, $n\ge m+2\ge 7$ and $2\nmid
mn$. Then $r(K_{1,m-1},T_n^1)=m+n-5.$
\endpro
 Proof. Since $n>m$ and $2\nmid mn$, we have $n\ge m+2$. Let $G$ be a graph of order $m+n-5$ such that $\overline{G}$ does
not contain any copies of $K_{1,m-1}$. We show that $G$ contains a
copy of $T_n^1$. Clearly $\Delta(\overline{G})\le m-2$ and so
$\delta(G)=m+n-6-\Delta(\overline{G})\ge n-4$.
 If $\Delta(G)=n-4$, then $G$ is a regular
graph of order $m+n-5$ with degree $n-4$ and so
$(m+n-5)(n-4)=2e(G).$ Since $m+n-5$ and $n-4$ are odd, we get a
contradiction. Thus, $\Delta(G)\ge n-3$.
 Assume that $v_0\in V(G)$,
$d(v_0)=\Delta(G)=n-3+c$, $\Gamma(v_0)=\{v_1,\ldots,v_{n-3+c}\}$,
$V_1=\{v_0\}\cup\Gamma(v_0)$ and
$V_1'=V(G)-V_1=\{u_1,u_2,\ldots,u_{m-3-c}\}$. Since $\delta(G)\ge
n-4$, for $v_i\in \Gamma(v_0)$ we have $1+|\Gamma(v_i)\cap
\Gamma(v_0)|+|V_1'|\ge d(v_i)\ge n-4$ and so $|\Gamma(v_i)\cap
\Gamma(v_0)|\ge n-5-(m-3-c)=n-m-2+c\ge c$.
\par We first assume that $|V_1'|=m-3-c\ge 2$. For
$i=1,2$ we have $|\Gamma(u_i)\cap\Gamma(v_0)|+|V_1'|-1\ge d(u_i)\ge
\delta(G)\ge n-4$ and so $|\Gamma(u_i)\cap\Gamma(v_0)|\ge
n-4+1-(m-3-c)=n-m+c\ge 2$. Hence $G$ contains a copy of $T_n^1$. If
$|V_1'|=1$, then $c=m-4\ge 1$. Since $d(u_1)\ge n-4>1$, we have
$u_1v_j\in E(G)$ for some $v_j\in\Gamma(v_0)$. Recall that
$|\Gamma(v_i)\cap \Gamma(v_0)|\ge c\ge 1$ for $v_i\in\Gamma(v_0)$.
$G$ must contain a copy of $T_n^1$. Now assume that $|V_1'|=0$. That
is, $c=m-3$ and $G=G[V_1]$. Since $d(v_0)=n-3+m-3\ge n-3+2$ and
$d(v_i)\ge n-4\ge 3$ for $v_i\in\Gamma(v_0)$, we see that
$G[\Gamma(v_0)]$ contains a copy of $2K_2$ and so $G$ contains a
copy of $T_n^1$.
\par By the above, $G$ contains a copy of $T_n^1$. Therefore
$r(K_{1,m-1},T_n^1)\le m+n-5$. From  Lemma 2.3,
$$r(K_{1,m-1},T_n^j)\ge
m-1+n-3-\f{1-(-1)^{(m-2)(n-4)}}2=m+n-5.$$ Hence
$r(K_{1,m-1},T_n^1)=m+n-5$ as claimed.
\par\q
\pro{Lemma 4.1} Let $m,n\in\Bbb N$, $n\ge 15$, $m\ge 7$, $n>m+1+\f
8{m-6}$ and $T_n\in\{T_n'',T_n''',T_n^3\}$. Let $G_m$ be a connected
graph of order $m$ such that $ex(m+n-5;G_m)\le \f{(m-2)(m+n-5)}2$.
Then $r(G_m,T_n)\le m+n-5.$ Moreover, if $m-1\mid n-5$, then
$r(G_m,T_n)=m+n-5$.
\endpro
Proof. If $T_n\not=T_n^3$ or $m\notin\{n-3,n-4\}$, appealing to
Lemmas 2.7 and 2.9 we have
$$\align\ex(m+n-5;T_n)
&=\f{(n-2)(m+n-5)-(m-4)(n-m+3)}2\\&\qq+
\max\Big\{0,\Big[\f{(m-4)(n-m)-3(n-1)}2\Big]\Big\}
\\&=\f{(n-2)(m+n-5)-(m-4)(n-m+3)}2\\&\qq+
\max\Big\{0,\Big[\f{(m-7)(n-m-3)-18}2\Big]\Big\}.
\endalign$$
 Thus, if $(m-7)(n-m-3)\ge 18$, then
 $$\align&\ex(m+n-5;G_m)+\ex(m+n-5;T_n)
\\&\le \f{(m-2)(m+n-5)}2+\f{(n-2)(m+n-5)-(m-4)(n-m+3)}2
\\&\qq+\f{(m-7)(n-m-3)-18}2
\\&=\f{(m+n-5)(m+n-7)}2<\b{m+n-5}2.
\endalign$$
 If $(m-7)(n-m-3)<18$,  since $n>m+1+\f 8{m-6}$ we see that
$(m-6)n>m^2-5m+2$, $(m-4)(n-m+3)>2(m+n-5)$ and so
$$\align&\ex(m+n-5;G_m)+\ex(m+n-5;T_n)
\\&\le\f{(m-2+n-2)(m+n-5)-(m-4)(n-m+3)}2
<\b{m+n-5}2.\endalign$$ Hence, $r(G_m,T_n)\le m+n-5$ by Lemma 2.1.
\par For $m=n-3$, using Lemma 2.8 we see that
$$\align &\ex(m+n-5;G_m)+\ex(m+n-5;T_n^3)
\\&=\ex(2n-8;G_{n-3})+\ex(2n-8;T_n^3)
 \\&\le \f{(2n-8)(n-5)}2 +n^2-9n+29+\max\big\{0,\big[\f{n-37}4\big]\big\}
\\&=2n^2-18n+49+\max\big\{0,\big[\f{n-37}4\big]\big\}
\\&<2n^2-17n+36=\b{m+n-5}2.
\endalign$$
 For $m=n-4$, appealing to Lemma 2.8,
$$\align &\ex(m+n-5;G_m)+\ex(m+n-5;T_n^3)
\\&=\ex(2n-9;G_{n-4})+\ex(2n-9;T_n^3)
 \\&\le \f{(2n-9)(n-6)}2 +n^2-10n+24
 +\max\Big\{\Big[\f n2\Big],13\Big\}
\\&=2n^2-20n+51-\f n2+\max\Big\{\Big[\f n2\Big],13\Big\}
\\&<2n^2-19n+45=\b{m+n-5}2.
\endalign$$
Thus,  $r(G_m,T_n^3)\le m+n-5$ for $m=n-4,n-3$ by Lemma 2.1.
\par Now assume that $m-1\mid n-5$. Then $m+n-6=k(m-1)$ for
$k\in\{2,3,\ldots\}$. Since $\Delta(\overline{kK_{m-1}})=n-5$ we see
that $kK_{m-1}$ does not contain $G_m$ as a subgraph and
$\overline{kK_{m-1}}$ does not contain $T_n$ as a subgraph. Hence
$r(G_m,T_n)>k(m-1)=m+n-6$ and so $r(G_m,T_n)=m+n-5$. The proof is
now complete.
\par\q
 \pro{Theorem 4.4} Let $m,n\in\Bbb N$, $n\ge 15$, $m\ge
7$, $n>m+1+\f 8{m-6}$ and $T_n\in\{T_n'',T_n''',T_n^3\}$. If $2\mid
m(n-1)$, then $r(K_{1,m-1},T_n)=m+n-5.$
\endpro
Proof. By Euler's theorem or Lemma 2.4, $\ex(m+n-5;K_{1,m-1})\le
\f{(m-2)(m+n-5)}2$. Thus, applying Lemma 4.1 we obtain
$r(K_{1,m-1},T_n)\le m+n-5$. Suppose that $2\mid m(n-1)$. By Lemma
2.3,
$$r(K_{1,m-1},T_n)\ge m-1+n-4-\f{1-(-1)^{(m-2)(n-5)}}2
=m+n-5.$$ Thus the result follows.
\par\q
 \pro{Corollary 4.5} Let $n\in\Bbb
N$, $n\ge 17$ and $T_n\in\{T_n'',T_n''',T_n^3\}$. Then
$r(K_{1,n-3},T_n)=2n-7$.
\endpro
Proof. Taking $m=n-2$ in Theorem 4.4 gives the result.
\par\q
\pro{Theorem 4.5} Let $m,n\in\Bbb N$, $m\ge 6$, $n\ge m+3$ and
$2\nmid m(n-1)$. Then
$$r(K_{1,m-1},T_n'')=r(K_{1,m-1},T_n''')=m+n-6.$$
\endpro
Proof. Let $G$ be a graph of order $m+n-6$ such that $\overline{G}$
does not contain any copies of $K_{1,m-1}$. That is,
$\Delta(\overline{G})\le m-2$. Thus,
$\delta(G)=m+n-7-\Delta(\overline{G})\ge n-5$.
 If $\Delta(G)=n-5$, then $G$ is a regular
graph of order $m+n-6$ with degree $n-5$ and so
$(m+n-6)(n-5)=2e(G).$ Since $m+n-6$ and $n-5$ are odd, we get a
contradiction. Thus, $\Delta(G)\ge n-4$.
 Assume that $v_0\in V(G)$,
$d(v_0)=\Delta(G)=n-4+c$, $\Gamma(v_0)=\{v_1,\ldots,v_{n-4+c}\}$,
$V_1=\{v_0\}\cup\Gamma(v_0)$ and
$V_1'=V(G)-V_1=\{u_1,\ldots,u_{m-3-c}\}$. Since $\delta(G)\ge n-5$,
we see that for
 $v_i\in \Gamma(v_0)$, $|\Gamma(v_i)\cap \Gamma(v_0)|+1+|V_1'|\ge
 d(v_i)\ge n-5$ and so
 $$|\Gamma(v_i)\cap \Gamma(v_0)|\ge
n-5-1-(m-3-c)=n-m-3+c\ge c.$$
 For $u_i\in V_1'$, we see that  $|\Gamma(u_i)\cap \Gamma(v_0)|+|V_1'|-1\ge
 d(u_i)\ge n-5$ and so
 $$|\Gamma(u_i)\cap \Gamma(v_0)|\ge
n-5-(m-4-c)=n-m-1+c\ge 2+c.$$
\par We first assume that $c=0$. Since $|V_1'|=m-3\ge 3$ and $\delta(G)\ge n-5$,
 we see that
$|\Gamma(u_i)\cap\{v_1,\ldots,v_{n-4}\}|\ge n-5-(m-4)=n-m-1\ge 2$
for $u_i\in V_1'$ and so $e(V_1V_1')\ge (m-3)(n-m-1)$. Since $n\ge
m+3$ we see that $(m-4)n\ge (m-4)(m+3)=m^2-m-12>m^2-2m-7$ and so
$e(V_1V_1')\ge (m-3)(n-m-1)>n-4$. Therefore, $|\Gamma(v_i)\cap
V_1'|\ge 2$ for some $i\in\{1,2,\ldots,n-4\}$. We may suppose that
$u_1v_i,u_2v_i,u_2v_j,u_3v_k\in E(G)$ with no loss of generality,
where $v_i,v_j,v_k$ are distinct vertices in $\Gamma(v_0)$. Thus $G$
contains a copy of $T_n^{''}$ and a copy of $T_n^{'''}$. \par Next
we assume that $|V_1'|=m-3-c\ge 3$ and $c\ge 1$. Then
$|\Gamma(u_i)\cap \Gamma(v_0)|\ge 3$ for $i=1,2,3$. Hence there are
distinct vertices $v_j,v_k,v_l\in\Gamma(v_0)$ such that
$u_1v_j,u_2v_k,u_3v_l\in E(G)$ and so $G$ contains a copy of
$T_n^{'''}$. Since $d(v_j)\ge n-5>2$, $v_j$ is adjacent to some
vertex $w$ different from $v_0$ and $u_1$. Hence, $G$ contains a
copy of $T_n^{''}$.
\par Now assume that $|V_1'|=2$. That is, $c=m-5$. Since
$|\Gamma(u_i)\cap \Gamma(v_0)|\ge \delta(G)-1\ge n-6\ge 3$ for
$i=1,2$, and $|\Gamma(v_i)\cap \Gamma(v_0)|\ge n-m-3+c=n-8\ge 1$ for
$v_i\in\Gamma(v_0)$, it is easy to see that $G$ contains a copy of
$T_n^{''}$ and a copy of $T_n^{'''}$.
\par Suppose that $|V_1'|=1$. Then $c=m-4\ge 2$, $d(u_1)\ge \delta(G)\ge
n-5\ge 4$ and $d(v_i)\ge \delta(G)\ge n-5\ge 4$ for
$i=1,2,\ldots,n-4+m-4$. Hence $G$ contains a copy of $T_n^{''}$ and
a copy of $T_n^{'''}$.
\par Finally we assume that $|V_1'|=0$. That is, $c=m-3$. Since
$d(v_i)\ge \delta(G)\ge n-5\ge 4$ for $i=1,2,\ldots,n-4+m-3$, it is
easy to see that $G$ contains a copy of $T_n^{''}$ and a copy of
$T_n^{'''}$.
\par Suppose that $T_n\in\{T_n^{''},T_n^{'''}\}$. By the above,
$G$ contains a copy of $T_n$. Hence
$r(K_{1,m-1},T_n)\le m+n-6$. By  Lemma 2.3, $r(K_{1,m-1},T_n)\ge
m-1+n-4-\f{1-(-1)^{(m-2)(n-5)}}2=m+n-6$. Thus
$r(K_{1,m-1},T_n)=m+n-6$ as asserted.
\par\q
\pro{Theorem 4.6} Let $n\in\Bbb N$ with $n\ge 15$. Then
$r(K_{1,n-4},T_n^3)=2n-8.$
\endpro
Proof. By Euler's theorem,
$\t{ex}(2n-8;K_{1,n-4})\le\f{(n-5)(2n-8)}2$. Thus,
$r(K_{1,n-4},T_n^3)\le 2n-8$ by taking $G_m=K_{1,n-4}$ in Lemma 4.1.
If $2\nmid n$, from Lemma 2.3 we have $r(K_{1,n-4},T_n^3)\ge
n-4+n-4=2n-8$. Thus the result is true for odd $n$. Now assume that
$2\mid n$. Let $G_0$ be the graph of order $2n-9$ constructed in
Theorem 3.2. Then $G_0$ does not contain $T_n^3$ as a subgraph. As
$\delta(G_0)=n-5$, we have $\Delta(\overline{G_0})=2n-10-(n-5)=n-5$
and so $\overline{G_0}$ does not contain $K_{1,n-4}$ as a subgraph.
Hence $r(K_{1,n-4},T_n^3)>|V(G_0)|=2n-9$ and so
$r(K_{1,n-4},T_n^3)=2n-8$ as claimed.
\par\q
 \pro{Theorem 4.7} Let $n\in\Bbb N$ with
$n\ge 10$. Then
$$r(K_{1,n-2},T_n^3)=r(K_{1,n-2},T_n'')=r(K_{1,n-2},T_n''')=2n-5.$$
\endpro
Proof. Let $T_n\in\{T_n'',T_n''',T_n^3\}$. Since
$\Delta(K_{1,n-2})=n-2$ and $\Delta(T_n)=n-4$,  we have
$r(K_{1,n-2},T_n)\ge 2(n-2)-1=2n-5$ by Lemma 2.3(ii). By Lemmas 2.4,
2.7 and 2.9,
$$\align &\ex(2n-5;K_{1,n-2})=\Big[\f{(n-3)(2n-5)}2\Big]=n^2-6n+8+
\Big[\f{n-1}2\Big],
\\& \ex(2n-5;T_n)=\f{(n-2)(2n-5)-3(n-4)}2=n^2-6n+11.
\endalign$$ Thus,
$$\align &\ex(2n-5;K_{1,n-2})+\ex(2n-5;T_n)
\\&=n^2-6n+8+\Big[\f{n-1}2\Big]+n^2-6n+11
<2n^2-11n+15=\b{2n-5}2.\endalign$$ Now, applying Lemma 2.1 yields
$r(K_{1,n-2},T_n)\le 2n-5$ and so $r(K_{1,n-2},$ $T_n)=2n-5$, which
proves the theorem.
\section*{5. Formulas for $r(T_m',T_n''),r(T_m',T_n''')$ \\and
$r(T_m',T_n^3)$}
\par\q
 \pro{Theorem 5.1} Let $m,n\in\Bbb N$, $n\ge 15$,
$m\ge 7$ and $m-1\mid n-5$. Suppose that
$G_m\in\{P_m,T_m',T_m^*,T_m^1,T_m^2,T_m^3,T_m'',T_m'''\}$ and
$T_n\in\{T_n'',T_n''',T_n^3\}$. Assume that $m\ge 10$ or $G_m
\not\in\{T_m^3,T_m'',T_m'''\}$. Then $r(G_m,T_n)=m+n-5.$
\endpro
Proof. Note that $m+n-5\e 1\mod {m-1}$. By (2.1) and Lemmas 2.5,
2.6, 2.10 and 2.11, $\ex(m+n-5;G_m)\le \f{(m-2)(m+n-5)}2$. Thus,
applying Lemma 4.1 and the fact $n\ge m+4$ gives the result.
\par\q
 \pro{Theorem 5.2} Let $m,n\in\Bbb N$, $m\ge 9$,
$n>m+2+\max\{0,\f {20-m}{m-8}\}$ and $m-1\nmid n-5$. Then
$$r(T_m',T_n'')=r(T_m',T_n''')=r(T_m',T_n^3)=m+n-6.$$
\endpro
Proof. Let $T_n\in\{T_n'',T_n''',T_n^3\}$. Since
$\Delta(T_m')=m-2<m-1$ and $\Delta(T_n)=n-4>m-2$,  we have
$r(T_m',T_n)\ge m-2+n-4=m+n-6$ by Lemma 2.3(ii)-(iii). Note that
$m\ge 9$ and so $n\ge 15$. Since $n>m+2+\f{20-m}{m-8}$, we see that
$(m-8)n>m^2-7m+4$ and so $(m-5)(n-m+4)>3(m+n-6)-(m-2)$.
 \par Suppose that
$T_n\not=T_n^3$ or $n\not=m+3$.
 From Lemmas 2.7 and 2.9, if $(m-5)(n-m+1)\ge 3(n-1)$, then
$$\align &\ex(m+n-6;T_n)
\\&\le \f{(n-2)(m+n-6)-(m-5)(n-m+4)}2
\\&\q+\f{(m-5)(n-m+1)-3(n-1)}2 =\f{(n-5)(m+n-6)}2;
\endalign$$
 if $(m-5)(n-m+1)<3(n-1)$, then
$$\align &\ex(m+n-6;T_n)\\&=\f{(n-2)(m+n-6)-(m-5)(n-m+4)}2
\\&<\f{(n-2)(m+n-6)-3(m+n-6)+m-2}2\\&=\f{(n-5)(m+n-6)+m-2}2.
\endalign$$
Recall that $m-1\nmid n-5$. By Lemma 2.5, $\ex(m+n-6;T_m')\le
\f{(m-2)(m+n-6)-(m-2)}2.$ Thus,
$$\align&\ex(m+n-6;T_m')+\ex(m+n-6;T_n)
\\&<\f{(m-2)(m+n-6)-(m-2)}2+\f{(n-5)(m+n-6)+m-2}2
\\&=\b{m+n-6}2.\endalign$$ Now applying Lemma 2.1 yields
$r(T_m',T_n)\le m+n-6$ and so $r(T_m',T_n)=m+n-6$.
\par Now assume that $T_n=T_n^3$ and $n=m+3$. Then
$\max\{0,\f {20-m}{m-8}\}<1$ and so $m=n-3\ge 15$. Also,
$m+n-6=2n-9=n-1+n-8=2m-3=m-1+m-2$. From Lemma 2.9(iii),
$$\ex(2n-9;T_n^3)=n^2-10n+24
 + \max\Big\{\Big[\f
n2\Big],13\Big\}.$$ By Lemma 2.5, $\ex(2m-3;T_m')
=\f{(m-2)(2m-3)-(m-2)}2=(m-2)^2=(n-5)^2$. Thus,
$$\align &\ex(m+n-6;T_m')+\ex(m+n-6;T_n^3)
\\&=(n-5)^2+n^2-10n+24
 + \max\Big\{\Big[\f
n2\Big],13\Big\}\\&=2n^2-20n+49+\max\Big\{\Big[\f n2\Big],13\Big\}
\\&<2n^2-19n+45=\b{2n-9}2.
\endalign$$
Applying Lemma 2.1 gives $r(T_m',T_n^3)\le m+n-6$ and so
$r(T_m',T_n^3)=m+n-6$ for $n=m+3$. This completes the proof.
\par\q

\pro{Theorem 5.3} Let $n\in \Bbb N$ with $n\ge 18$. Then
$$r(T_{n-3}',T_n'')=r(T_{n-3}',T_n''')=r(T_{n-3}',T_n^3)=2n-9.$$
\endpro
Proof. Suppose that $T_n\in\{T_n'',T_n''',T_n^3)$. Since
$\Delta(T_n)=n-4>n-5=\Delta(T_{n-3}')$, from Lemma 2.3(ii) we have
$r(T_{n-3}',T_n)\ge 2(n-4)-1=2n-9$. By Lemma 2.5,
$\ex(2n-9;T_{n-3}') =\f{(n-5)(2n-10)}2=n^2-10n+25$. From Lemma 2.7,
for  $T_n\in\{T_n'',T_n'''\}$,
$$\align &\ex(2n-9;T_n)\\&=\f{(n-2)(2n-9)-7(n-8)}2
+\max\Big\{0,\Big[\f{4(n-8)-3(n-1)}2 \Big]\Big\}
\\&=n^2-10n+37+\t{max}\Big\{0,\Big[\f{n-29}2\Big]\Big\}<n^2-9n+20\endalign$$
and so
$$\ex(2n-9;T_{n-3}')+\ex(2n-9;T_n)
<n^2-10n+25+n^2-9n+20 =\b{2n-9}2.$$ Now, applying Lemma 2.1 yields
$r(T_{n-3}',T_n)\le 2n-9$ and so $r(T_{n-3}',T_n)$ $=2n-9$. On the
other hand, from Lemma 2.8 we have
$$\ex(2n-9;T_n^3)=n^2-10n+24+\max\Big\{\Big[\f n2\Big],13\Big\} <n^2-9n+20.$$
 Thus, $$\ex(2n-9;T_{n-3}')+\ex(2n-9;T_n^3)<n^2-10n+25+n^2-9n+20
 =\b{2n-9}2.$$
Applying Lemma 2.1, $r(T_{n-3}',T_n^3)\le 2n-9$ and so
$r(T_{n-3}',T_n^3)=2n-9$, which completes the proof.
\par\q
 \pro{Theorem 5.4} Let $m,n\in \Bbb N$ with $n>m\ge 10$, and
$T_m\in\{T_m'',T_m''',$ $T_m^3\}$. Then
$$\align r(T_m,T_n')&=r(T_m,T_n^*)\\&=\cases m+n-3&\t{if $m-1\mid n-3$,}
\\m+n-4&\t{if $m-1\nmid n-3$ and $n\ge (m-3)^2+2$.}
\endcases\endalign$$
\endpro
Proof. If $m-1\mid n-3$, then
$\ex(m+n-3;T_m)=\f{(m-2)(m+n-3)-(m-2)}2$ by Lemmas 2.7 and 2.9.
Thus, the result follows from [9, Theorems 4.1 and 5.1].
\par Now assume that $m-1\nmid n-3$. By Lemma 2.10, $\ex(m+n-4;T_m)
<\f{(m-2)(m+n-4)}2$. Applying [9, Theorems 4.4 and 5.4] deduces the
result. The proof is now complete.

\section*{6. Evaluation of $r(T_m^0,T_n)$ with
$T_m^0\in\{T_m^*,T_m^1,$ \\ $T_m^2\}$ and
$T_n\in\{T_n'',T_n''',T_n^3\}$}
\par\q
 \pro{Lemma 6.1 ([7, Theorem 8.3, pp.11-12])}
Let $a,b,n\in\Bbb N$. If $a$ is coprime to $b$ and $n\ge
(a-1)(b-1)$, then there are two nonnegative integers $x$ and $y$
such that $n=ax+by$.\endpro
\par\q
 \pro{Theorem 6.1} Let $m,n\in \Bbb N$ with
 $m\ge 9$,  $n> m+1+\f{12}{m-8}$ and $m-1\nmid n-5$.
 Suppose that $T_m^0\in\{T_m^*,T_m^1,T_m^2\}$ and
$T_n\in\{T_n'',T_n''',T_n^3\}$. Assume that $T_m^0\not=T_m^*$ or
$m\ge 11$. Then $r(T_m^0,T_n)=m+n-7$ or $m+n-6$. If $n\ge (m-3)^2+4$
or $m+n-7=(m-1)x+(m-2)y$ for some nonnegative integers $x$ and $y$,
then $r(T_m^0,T_n)=m+n-6$.
\endpro
Proof. Note that $\Delta(T_m^0)=m-3<n-4=\Delta(T_n)$. Using Lemma
2.3(ii)-(iii), $r(T_m^0,T_n)\ge m-3+n-4=m+n-7$.  Since $m-1\nmid
n-5$, from Lemmas 2.6, 2.11 and 2.12 we have $\ex(m+n-6;T_m^0)\le
\f{(m-2)(m+n-6)-(m-2)}2$.
\par We first assume that $T_n\not=T_n^3$ or $n\not=m+2,m+3$. By the
proof  of Theorem 5.2, $\ex(m+n-6;T_n)<\f{(n-5)(m+n-6)+m-2}2$. Thus,
$$\align &\ex(m+n-6;T_m^0)+\ex(m+n-6;T_n)
\\&<\f{(m-2)(m+n-6)-(m-2)}2 +\f{(n-5)(m+n-6)+m-2}2\\&=\b{m+n-6}2.\endalign$$
Hence, $r(T_m^0,T_n)\le m+n-6$ by Lemma 2.1 and so
$r(T_m^0,T_n)=m+n-6$ or $m+n-7$.
\par We next assume that $T_n=T_n^3$ and $n=m+2$. Then $m+n-6=2n-8=n-1+n-7$,
$m+2>m+1+\f{12}{m-8}$ and so $n-2=m>20$. By Lemma 2.9(iv),
$$\align &\ex(m+n-6;T_n^3)\\&=\ex(2n-8;T_n^3)=\f{(n-2)(2n-8)-6(n-7)}2
+\max\Big\{\Big[\f{n-37}4\Big],0\Big\}
\\&=n^2-9n+29+\max\Big\{\Big[\f{n-37}4\Big],0\Big\}<n^2-9n+29+\f{n-22}2.
\endalign$$ Thus,
$$\align &\ex(m+n-6;T_m^0)+\ex(m+n-6;T_n)
\\&< \f{(n-4)(2n-9)}2+n^2-9n+29+\f{n-22}2\\&=(n-4)(2n-9)=\b{2n-8}2.
\endalign$$
Hence $r(T_m^0,T_n^3)\le m+n-6$ by Lemma 2.1 and so
$r(T_m^0,T_n^3)=m+n-6$ or $m+n-7$.
\par Finally, we assume that $T_n=T_n^3$ and $n=m+3$.
Then $m+n-6=2n-9=n-1+n-8$, $m+3>m+1+\f{12}{m-8}$ and so $n-3=m\ge
15$. From Lemma 2.9(iii),
$$\align &\ex(m+n-6;T_n^3)\\&=\ex(2n-9;T_n^3)=\f{(n-2)(2n-9)-7n+30}2
+\max\Big\{\Big[\f{n}2\Big],13\Big\}
\\&=n^2-10n+24
 + \max\Big\{\Big[\f
n2\Big],13\Big\}.\endalign$$  Recall that
$$\align \ex(m+n-6;T_m^0)=
\ex(2m-3;T_m^0)
&\le\f{(m-2)(2m-3)-(m-2)}2\\&\q=(m-2)^2=(n-5)^2.\endalign$$ We then
obtain
$$\align &\ex(m+n-6;T_m^0)+\ex(m+n-6;T_n^3)
\\&=(n-5)^2+n^2-10n+24
 + \max\Big\{\Big[\f
n2\Big],13\Big\}\\&=2n^2-20n+49+\max\Big\{\Big[\f n2\Big],13\Big\}
\\&<2n^2-19n+45=\b{2n-9}2.
\endalign$$
Applying Lemma 2.1 gives $r(T_m^0,T_n^3)\le m+n-6$ and so
$r(T_m^0,T_n^3)=m+n-6$ or $m+n-7$ for $n=m+3$.
\par If $m+n-7=(m-1)x+(m-2)y$ for some
nonnegative integers $x$ and $y$, setting $G=xK_{m-1}\cup yK_{m-2}$
we find that $G$ does not contain any copies of $T_m^0$. Observe
that $\Delta(\overline{G})=n-5$ or $n-6$. We see that $\overline{G}$
does not contain any copies of $T_n$. Hence
$r(T_m^0,T_n)>|V(G)|=m+n-7$ and so $r(T_m^0,T_n)=m+n-6$. If $n\ge
(m-3)^2+4$, then $m+n-7\ge (m-2)(m-3)$. By Lemma 6.1,
$m+n-7=(m-1)x+(m-2)y$ for some nonnegative integers $x$ and $y$ and
so $r(T_m^0,T_n)=m+n-6$ as claimed.
\par Summarizing the above proves the theorem.

\end{document}